\documentclass[10pt,onecolumn]{IEEEtran} 

\usepackage{cite}

\ifCLASSINFOpdf
\usepackage[pdftex]{graphicx}
\else
\fi
\usepackage[cmex10]{amsmath}
\usepackage{amsfonts}
\usepackage{amssymb}

\usepackage{mdwmath}
\usepackage{mdwtab}
\usepackage{booktabs}
\usepackage[font=footnotesize]{caption}
\usepackage{subcaption}

\usepackage{stfloats}

\usepackage{pgfplots}
\usepackage{mathtools}

\begin{document}

\hyphenation{op-tical net-works semi-conduc-tor}

\title{Robust Corrective Control Measures in Power Systems with Dynamic Line Rating}

\author{\IEEEauthorblockN{Matthias A. Bucher and G$\ddot{\textnormal{o}}$ran Andersson\\}
\IEEEauthorblockA{Power Systems Laboratory, ETH Zurich \\
Physikstrasse 3, 8092 Zurich\\ $\left\{ \textnormal{mb, andersson} \right\}$@eeh.ee.ethz.ch}}

\maketitle

\begin{abstract}
Dynamic Thermal Line Rating (DLR) is deemed to be an effective way to increase transmission capacities and therefore enabling additional operational flexibility.
The transmission capacities are dynamically determined based on current or expected weather conditions. First pilot projects have proven its efficiency.

In this paper we present two approaches to determine the location and amount of ramping capabilities for corrective control measures in the case of errors in the forecast of the line rating.
Both approaches result in robust optimization problems, where the first approach guarantees that there is a suitable remedial action for every realization of uncertainty in a given uncertainty set. The corrective control action is calculated once the forecast error is known. The second approach relies on affine policies which directly relate the current line rating to corrective control measures, if needed, which enables a decentralized operation.
In case studies, the two approaches are compared and the reduction of overall operational costs is demonstrated exemplarily for different parameters.
\end{abstract}


\IEEEpeerreviewmaketitle

\setlength{\tabcolsep}{1.1em}

\section{Introduction}

Many power systems are currently undergoing fundamental changes. For example in Europe, the shares of fluctuating renewable energy sources are increasing for a sustainable electricity production and regional markets are merged together or coupled \cite{IEMDirective}. Both developments result in power flows over longer distances and power flow patterns changing more frequently. In order to guarantee a secure and reliable operation under these conditions, on one hand, the transmission system needs to be operated in a more flexible way and substantial increases in transmission capacities are essential. Network reinforcement is a possible remedy. The drawbacks are long licensing procedures and public opposition as well as high investment costs. Therefore also alternatives should be considered, i.e. methods are sought, that allow to use the full potential of transmission lines. One approach, Dynamic Line Rating (DLR), is to set the line ratings according to current meteorological conditions. Traditionally, the nominal line ratings (NLR) are based on different operational constraints, in the transmission system mainly on the thermal limitations for given worst-case weather conditions. By monitoring the thermal state of the transmission lines and considering meteorological conditions, dynamically adapted line ratings enable additional transmission capacity. One major challenge is to handle the uncertainty in the forecast of the ampacity.

DLR has been investigated in different pilot projects, e.g. \cite{Sattinger}, and is considered by many grid operators in the development plans of their grids. In \cite{DLRBolun} the potential of dynamic line ratings with respect to the integration of fluctuating renewable energy sources is investigated. Other authors have shown how to operate DLR in combination with connections to wind farms \cite{McLaughlin2011,windpark,windcooling}. In previous work \cite{MB2}, we have shown how dynamic line ratings can be incorporated in a probabilistically N-1 secure dispatch.


In this paper we present two approaches that allow to allocate additional transmission capacity based on meteorological forecasts. The methods trade potentially higher line ratings based on current meteorological forecasts off against the costs arising from remedial actions when transmission lines would be overloaded due to wrongly forecasted line ratings. The methods optimize the procurement in terms of costs and consider the location of injection in the grid. The capacities for re-dispatching could be provided by flexible units, e.g. such as fast-ramping hydro-storage power plants, or from demand side participation. Also thinkable would be the provision by devices, that allow to control the power flows on a certain line, e.g. HVDC lines. However, this is not within the scope of this paper.

The approaches are based on robust optimization techniques and guarantee secure operation for all realizations of possible line ratings within a given uncertainty set.
A first approach relies on a central re-dispatch coordination entity, e.g. the TSO, for the operation of the DLR system, which leads to lower operational costs compared to the second approach but also needs central computation possibilities and reliable communication infrastructure to the entities providing the reserves. The second approach is based on affine policies, i.e. the providers of reserves operate according to a policy relating the current line rating to a predefined remedial action. This approach can be operated in a decentralized manner. Both methods are compared in a case study based on the IEEE RTS96 system. The amount of procured reserves and the changes in cost for operating the system for different combinations of parameters are determined and their advantages and drawbacks summarized.

The paper is structured as follows: In section II, the modeling of the thermal behavior of the lines is shown and the potential and risk with respect to the ampacity distributions are discussed.
Section III introduces the general problem and IV presents the modeling of the uncertainty sets. In sections V and VI the methods are presented and finally applied to a case study in section VII.
The paper concludes in section VIII.


\section{Modeling, Potential and Risks of DLR}
\subsection{Modeling}
For this paper we use the IEEE model \cite{DLRIEEE} to estimate the temperature of the conductor. Although different models are discussed in literature, e.g. \cite{DLRcigre,DLRIEC}, most of them are based on the balance of cooling and heating effects affecting the line's temperature. This balance can be written as $q_c + q_r = q_s + I^2 \cdot R$. The term $q_c$ is the convection heat loss rate, depending on the wind speed and the angle between the wind direction and the line direction. $q_r$ is the radiated heat loss rate, which is a function of the difference between the conductor temperature and the ambient temperature. $q_s$ is the rate of solar heat gain caused by solar irradiation. The heating from the power flow depends on the current $I$  and the temperature dependent resistance of the conductor $R$. The steady-state equation is independent of the line length, i.e. the meteorological conditions are assumed to be uniform along the line.

\begin{figure}[!h]
\centering
\includegraphics[width=3.5in]{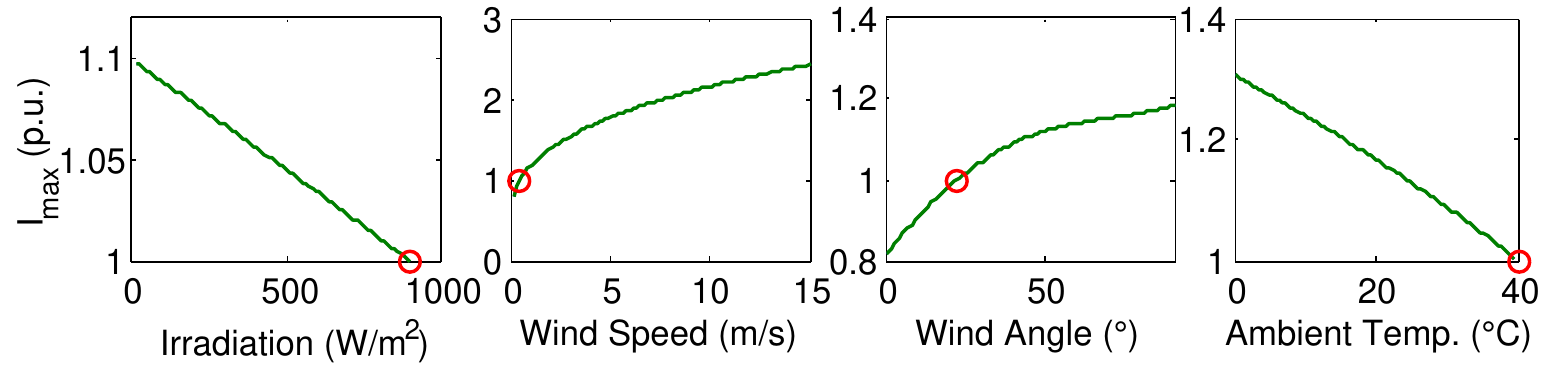}
\caption{Sensitivity of the ampacity of a "Drake" 26/7 ACSR line \cite{DLRIEEE} with respect to different meteorological conditions. Marked with red circle are the assumptions made for the calculation of the nominal line rating (NRL) such that the probability of an overloading is less than 0.5\% in the case of Switzerland. The nominal line rating is 1 p.u.}
\label{Fig:infl}
\end{figure}

Fig. \ref{Fig:infl} shows the influences of meteorological quantities on the ampacity by varying only one variable from the base case. As seen, the influence is substantial, especially the additional cooling due to higher wind speeds. In some cases, the ampacity is more than twice as high as in the base case.

The variables $q_r$, $q_s$ and $q_c$ depend on the current meteorological conditions and therefore also the ampacity $I_{max}$ of the conductor. The ampacity can be calculated for given meteorological conditions by selecting the maximum conductor temperature $\bar{T}$ and solving the balance equation for $I$:

\begin{equation}
\begin{aligned}
I_{max} \left( \bar{T} \right) = \sqrt{ \frac{q_c \left( \bar{T} \right) + q_r \left( \bar{T} \right)  - q_s }{R \left( \bar{T} \right)}} \\
\end{aligned}
\label{eq:dlr_steadystate_curr}
\end{equation}

The $\bar{T}$ is selected based on the thermal limits of the conductor, the allowed sag of the line as well as the degradation. The line rating can be calculated based on the ampacity and the voltage level of the transmission line.

\subsection{Potential and Risk}

As the static limits are calculated under conservative meteorological assumptions, the effective transmission capacity as function of the actual meteorological conditions will be higher in most of the cases.

\begin{figure}
\centering
\begin{minipage}{.4\textwidth}
  \centering
  \includegraphics[width=.6\linewidth]{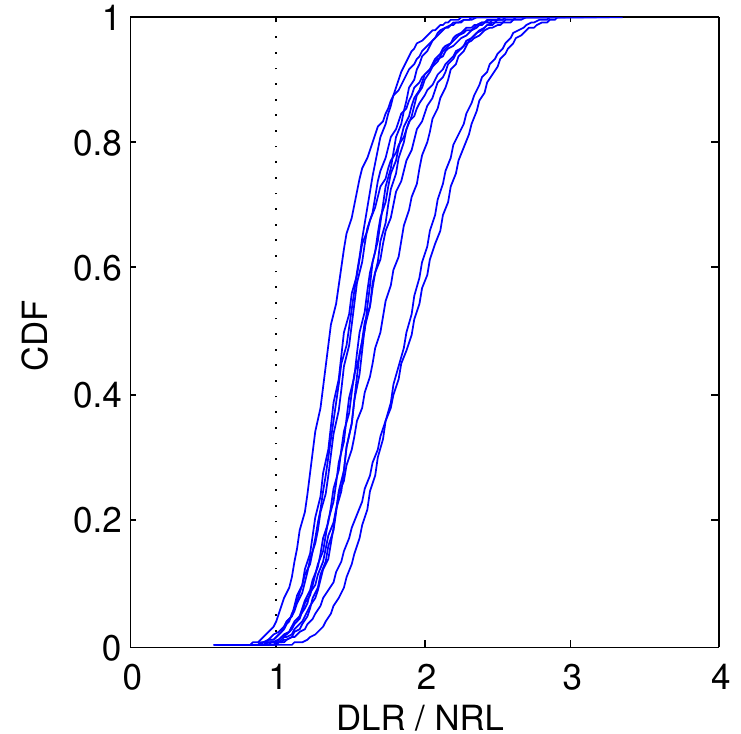}
  \captionof{figure}{Cumulative Distribution Function (CDF) of normalized line ratings based on meteorological data from different stations.}
  \label{fig:cdf}
\end{minipage}%
\hspace{1cm}
\begin{minipage}{.4\textwidth}
  \centering
  \includegraphics[width=.6\linewidth]{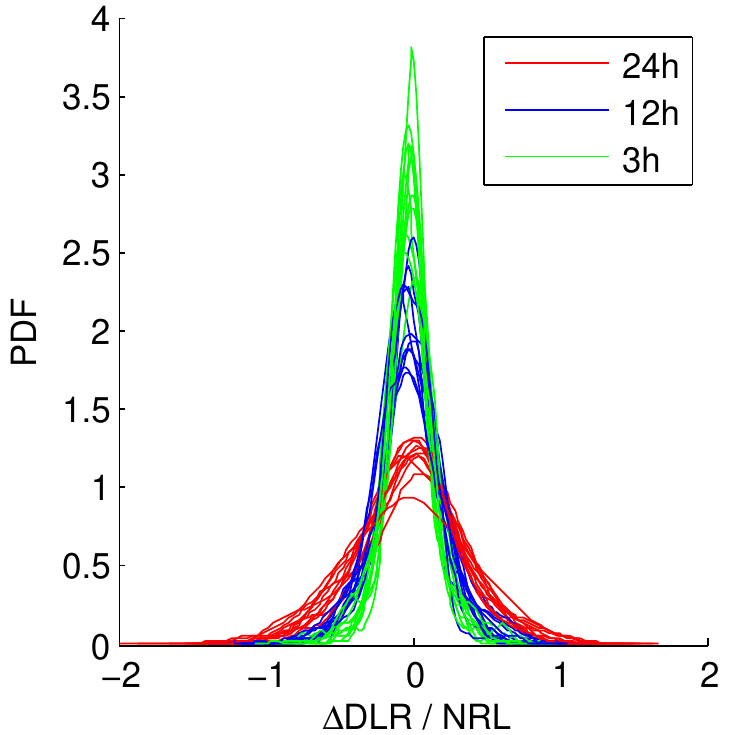}
  \captionof{figure}{Probability Density Function (PDF) of ampacity forecast errors for meteorological data from different stations and different forecast lead times normalized to nominal line ratings.}
  \label{fig:pdist}
\end{minipage}
\end{figure}

Fig. \ref{fig:cdf} displays the ascending theoretical ampacities calculated using realistic meteorological conditions \footnote[1]{Data from 11 weather stations in Switzerland from August 2011 until August 2012. The data contain measurements of wind, irradiation and ambient temperature every 3h and forecasts for the next 24h starting every 3h. Wind angle is not measured. We assume $22.5^{\circ}$ for our calculations.}. The values are normalized to the ampacitiy of the base case (NLR). It can be observed that the potential is significant. It should however be noted, that meteorological conditions are local quantities and the ampacity corresponds to the worst case along the whole line. Thus, for a successful implementation of DLR, a reasonable number of measurement stations would have to be placed along the line. We assume that the given meteorological data is representative for the whole line.

%
%

As the dispatch planning is typically performed day-ahead, only the meteorological forecasts are available.
The forecasts are not perfect, thus the ampacity forecast are uncertain. In Fig. \ref{fig:pdist} the distribution of the ampacity error, i.e. the difference of the ampacity calculation with forecasted and actual data\footnotemark[1] for different forecast horizons is displayed. It can be seen that the forecasted ampacity can differ quite significantly from the actual but the variation reduces for smaller forecast periods.

Based on the current we can calculate the thermal power limits of the lines for given meteorological conditions. The line ratings, denoted by the vector $\delta$, are random variables approximated by a multivariate normal distribution $\delta \thicksim \mathcal{N}\left(\mu, \Sigma \right)$, where $\mu$ is a vector with the expected line ratings and $\Sigma$ is the corresponding covariance. Possible correlations between different line ratings can be considered in the off-diagonal elements.

\section{General Problem Description}

\begin{figure}[!h]
\centering
\includegraphics[width=6.5cm]{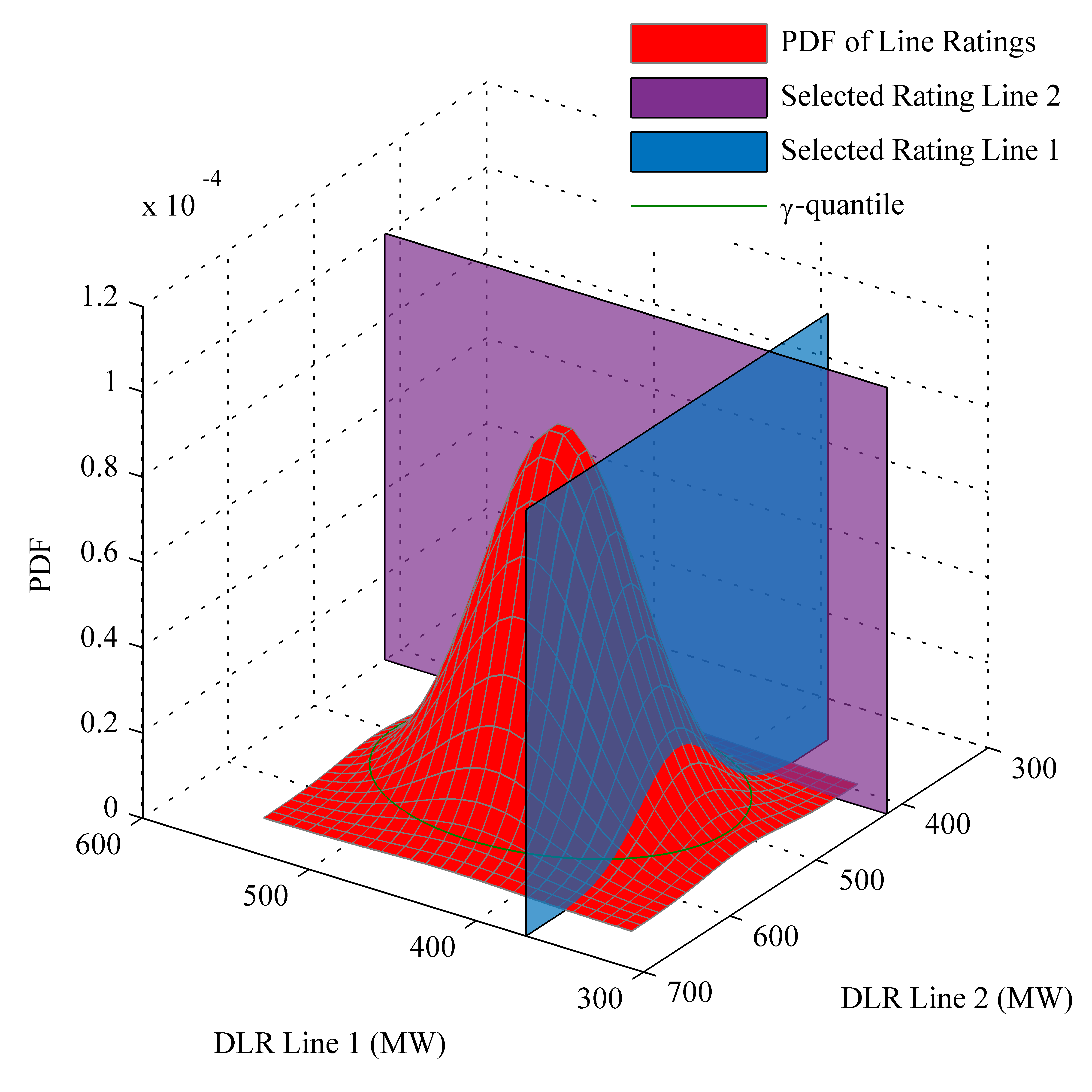}
\caption{Probability density function (PDF) of line ratings for two lines and $\gamma = 95\%$.}
\label{fig:illustuc}
\end{figure}

Fig. \ref{fig:illustuc} illustrates the probability distribution function (PDF, in red) for two lines. This PDF is based on the forecast for the considered period of time and might look different in other times. The cuts in blue represent selected line ratings providing upper operational limits for the corresponding power flows. For line 1, the selected line rating would be higher than the actual line rating $\delta$ for some realizations of $\delta$ and thus could be overloaded in some cases and would require alleviating actions.

Therefore, the line ratings have traditionally been selected such that most of the scenarios in all time periods are larger than the NLR. The methods presented in this paper allow to trade off the benefits of higher line ratings and the costs for remedial actions that might be necessary depending on the realization of $\delta$. As the realization of $\delta$ is unknown at the instant of planning the operation and can only be forecasted, sufficient capacity to avoid overloading for all possible realizations has to be procured. Mathematically, this can be written as the following constraint:

\begin{equation}
\begin{aligned}
& \left| P_L + f \left( \Delta, \delta \right) \right| \leq \delta \ \forall \delta \in \mathcal{W}\\
& \Delta \in \mathcal{R} \\
\end{aligned}
\label{genform}
\end{equation}

$P_L$ is the scheduled power flow and $f \left( \Delta, \delta \right)$ is a function that describes the change of power flow caused by a suitable remedial action $\Delta$ that is contained by the set $\mathcal{R}$ with all feasible remedial actions. The line ratings $\delta$ can take any value within the uncertainty set $\mathcal{W}$.
The goal of the methods is to procure sufficient reserves, i.e. solve a robust optimization problem, such that the set $\mathcal{R}$ is large enough to provide a suitable remedial action for every realization of uncertainty.
In the next section, the goal is to characterize the uncertainty sets $\mathcal{W}$ containing all relevant realizations of line ratings.

\section{Modeling of Uncertainty in Uncertainty Sets}
We now formulate uncertainty sets which define the bounds on possible realizations of line ratings $\delta$. This uncertainty set is dependent on different parameters, such as the forecast and the forecast uncertainty. For simplicity, we assume that the distribution for a given forecast is given as  $\delta \thicksim \mathcal{N}\left(\mu, \Sigma \right)$.
Considering again Fig. \ref{fig:illustuc}, the green ellipsoid bounds the $\gamma$-quantile, i.e. $\gamma$\% of all realizations are within the ellipsoid; in this case 95\%. The goal of this section is to define the sets containing all realization within the $\gamma$-quantile.

We are using two common types of uncertainty sets, {\it ellipsoidal} and {\it polyhedral} uncertainty \cite{BenTal99}. The polyhedral uncertainty will be used in approach I and ellipsoidal uncertainty will be used in approach II.

\subsubsection{Ellipsoidal Uncertainty}
We use ellipsoidal uncertainty sets due to their close relation to the normal distribution. The covariance matrix can be decomposed using the singular value decomposition $\Sigma = U \kappa U^T $ into: 

\begin{equation}
\Sigma  = U \kappa U^T = B B^T \text{ where } B = U \kappa^{\frac{1}{2}} \text{ and $\kappa$ is diagonal.}\\
\end{equation}

Introducing a random variable with the same dimensions as $\delta$ that is distributed standard normal, $z \thicksim \mathcal{N} \left( 0, I \right)$, the following holds:

\begin{equation}
\begin{aligned}
B^{-1}\left(\delta - \mu \right) \thicksim \mathcal{N} \left( 0, B^{-1}\Sigma \left( B^{-1} \right)^T \right) \Rightarrow \delta = B z + \mu\\
\end{aligned}
\end{equation}

Therefore the ellipsoidal uncertainty set is written as: $\mathcal{W} = \left\{ z: \mu + B z, \lVert z\rVert_2 \leq \sqrt{\rho} \right\}$.
As $z^T z = \left(\delta - \mu \right)^T \Sigma^{-1} \left(\delta - \mu \right) \thicksim \chi^2_k$, where $k$ is the dimension of $\delta$, the parameter $\rho$ is selected to be the $\gamma$-quantile of the $\chi^2_k$-distribution.

\subsubsection{Polyhedral Uncertainty}

Ellipsoidal uncertainty sets complicate robust optimization problems (see Eq. \eqref{eq:opt2}) and also, the distribution might not be normal as assumed in this paper. Therefore as an alternative, polyhedral uncertainty sets can be used, which can also be created based upon a selected number of scenarios, e.g. when the distribution is unknown.

\begin{figure}[!h]
\centering
\includegraphics[width=75pt]{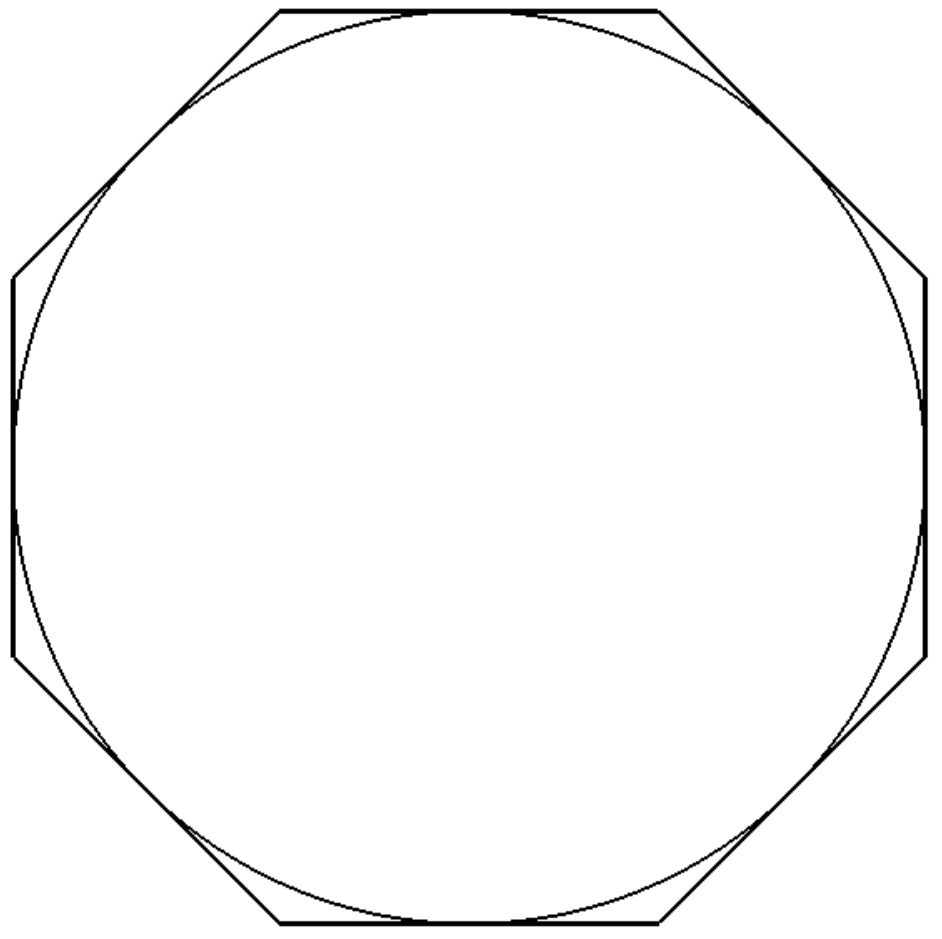}
\caption{Conservative approximation of ball-constraints using halfspaces.}
\label{fig:circleapprox}
\end{figure}

Based on the formulation above, we can approximate the ellipsoidal uncertainty by approximating the norm-constraint of $\mathcal{W}$ with halfspaces, i.e. we write $\mathcal{W}_P = \left\{ z: \mu + B z, S z \leq h \right\}$. The matrix inequality $S z \leq h$ represents the approximation of the ball-constraint in $\mathcal{W}$ by halfspaces as illustrated in Fig. \ref{fig:circleapprox}.

The resulting uncertainty set $\mathcal{W}_P$ is linear, which reduces complexity of the optimization problem.

\section{Approach I: Corrective Control using online optimization}
In this first approach, we assume that the balancing responsible party, e.g. the TSO, needs to procure a certain amount capacity from flexible units that can be instructed accordingly to alleviate congested lines. The goal of this approach is to find the suitable amounts and locations of such flexibility, that for any given uncertainty, the TSO is guaranteed to have sufficient capacity available for a suitable re-dispatching action. The procurement happens during the planning phase, however the action to take is only calculated during operation when the effective line rating is measured. As the TSO has to compute the corresponding action, this approach is a centralized one, which necessitates certain communication and computing infrastructure.

\subsection{Procurement and Dispatch}
The procurement of reserves is shown here in combination with a dispatch problem, but the procurement can also be performed independent of the dispatch problem.
In this case, we can apply a vertex enumeration technique, i.e. we procure sufficient flexibility, that every realization of $\delta$ which corresponds to a vertex of the polyhedral uncertainty set $\mathcal{W}_P$ is feasible. Due to convexity of $\mathcal{W}_P$ and the constraint formulation, every realization within the set is feasible as well. The amount of capacity to procure corresponds to the maximum, respectively minimum, of the needed capacities of all tested vertices.
This method has the advantage that it is simple and can be parallelized and directly incorporated in a dispatch algorithm. Drawbacks are that the expected costs of operation are not directly incorporated in the optimization problem, although expectation of certain realizations can be assumed to be known. Although that a action to keep the system secure is available, this action still needs to be calculated and applied during operation.
We now formulate the combined procurement and dispatch problem:

\newcommand*\circled[1]{\tikz[baseline=(char.base)]{
            \node[shape=circle,draw,inner sep=1.5pt] (char) {#1};}}

\begin{equation}
\begin{aligned}
c_{DLR}^* =& \min_{P_{gen},\underline{\Delta},\overline{\Delta},\Delta^-_i,\Delta^+_i} \  c_{disp} + c_{proc} + \alpha c_{wc}\\
&s.t. \\
\circled{A}&\left\{
\begin{aligned}
& c_{disp} \geq P_{gen}^T Q_{gen} P_{gen} + c_{gen}^T P_{gen} \\
& \sum_j \left( P_{gen,j} - P_{load,j}\right) = 0 \\
& \underline{P}_{gen} \leq P_{gen} \leq \overline{P}_{gen}
\end{aligned}
\right.\\
& \text{and } \forall\ i \text{, where $\delta_i$ is $i$-th vertex of } \mathcal{W}_P \\
\circled{B}&\left\{
\begin{aligned}
& c_{proc} \geq c^{+} \overline{\Delta} - c^{-} \underline{\Delta} \\
& c_{wc} \geq c_{op}^{+} \Delta^+_i - c_{op}^{-}\Delta^-_i  \\
&\underline{\Delta}^{max} \leq \underline{\Delta} \leq \Delta^-_i \leq 0 \\
& 0 \leq \Delta^+_i \leq \overline{\Delta} \leq \overline{\Delta}^{max}\\
& \sum  \Delta^+_i + \Delta^-_i = 0\\
\end{aligned}
\right.\\
\circled{C}&\left\{
\begin{aligned}
& \left| H_{DLR} \left(P_{gen} - P_{load}+ \Delta^+_i + \Delta^-_i \right) \right| \leq \delta_i \\
& \left| H_{NLR} \left(P_{gen} - P_{load}+ \Delta^+_i + \Delta^-_i \right) \right| \leq \overline{P}_{L,NLR} \\
\end{aligned}
\right.\\
\end{aligned}
\label{procI}
\end{equation}

The optimization minimizes the costs arising from the dispatch, the procurement of reserves and the worst-case operational costs of the reserves. The worst-case costs are weighted with the tuning parameter $\alpha$. The dispatch \circled{A} is constrained to the physical limits of the generation units $P_{gen}$ and the power balance. The costs are composed by a quadratic cost matrix $Q_{gen}$ and a linear cost vector $c_{gen}$. For every vertex of the uncertainty set, a set of constraints for the operation of the reserves is introduced in \circled{B}. The variables $\Delta^+_i,\ \Delta^-_i$ representing the remedial actions, that should be within the procured capacity $\underline{\Delta},\ \overline{\Delta}$, are selected such that for the given realization $\delta_i$ the lines are not overloaded (constraint \circled{C}) and the power balance is still satisfied. The amount of reserves to procure $\underline{\Delta},\ \overline{\Delta}$ is then given by the maximal activated amount of every possible provider of reserves, but limited to the maximum amounts that can be procured $\underline{\Delta}^{max},\overline{\Delta}^{max}$. The costs are calculated using the cost vectors $c^{+}, c^{-}$ for the procurement and $c_{op}^{+},c_{op}^{-}$ for the operation. All cost vectors are selected positive, as a redispatch always results in additional costs. The power flows are calculated using the PTDF matrix $H$, which is the sensitivity of the power flows with respect to the net bus injections \cite{Wollenberg2000}. $H$ is split into two matrices, $H_{DLR},H_{NLR}$, which correspond to influences of net bus injections on lines with DLR and without DLR, i.e. NLR. The corresponding power flow limits are given by $\delta$ and $\overline{P}_{L,NLR}$ respectively.

In general, this flexibility can be provided by flexible generation units, but also by demand side participation or technologies that allow to control the power flow, e.g. FACTS devices or HVDC interconnections. Without loss of generality, we assume that the ramping capacities are provided only by conventional generation units. However, the framework can be extended to account for other means as well.

\subsection{Operation of DLR}
During the operation, the TSO supervising the lines with DLR would have to calculate the remedial action once the deviation from the forecast is known. It would therefore solve the following problem to minimize the operational costs but avoiding overloaded lines.

\begin{equation}
\begin{aligned}
& \min_{\Delta^-, \Delta^+}\  c_{op}^{+} \Delta^+ - c_{op}^{-}\Delta^-\\
&s.t.\ \delta \text{ is {\it current} transmission capacity} \\
&\left\{
\begin{aligned}
&0 \leq \Delta^+  \leq \overline{\Delta}\\
&\underline{\Delta} \leq \Delta^- \leq 0\\
& \sum  \Delta^+ + \Delta^- = 0\\
& \left|  H_{DLR} \left(P_{gen} - P_{load}+ \Delta^+_i + \Delta^-_i  \right) \right| \leq \delta \\
& \left| H_{NLR} \left(P_{gen} - P_{load}+ \Delta^+_i + \Delta^-_i \right) \right| \leq \overline{P}_{L,NLR} \\
\end{aligned}
\right.
\end{aligned}
\end{equation}

$P_{gen},\ \overline{\Delta},\ \underline{\Delta}$ are not optimization variables as they are fixed according to the solution to \eqref{procI}. The remedial actions $\Delta^+, \Delta^-$ are communicated to the providers of reserves.

\section{Approach II: Corrective Control using affine policies}

In the second approach we assume an affine policy for the operation of the flexible units, i.e. the changes of generator setpoints are affine functions of the uncertainty realizations of the transmission capacities $\delta$. The concept of affine policies is widely used in reserves operations, e.g. the operation of primary control reserves is an affine function of the frequency deviation and units providing secondary control reserve rely on a signal of the TSO, with which they linearly scale their offered reserves. 
In this paper, the policies activate corrective measures given that the line capacity falls below a certain threshold. The units' change in active power output $\Delta$ is set to:

\begin{equation}
\begin{aligned}
& \Delta = \left\{
\begin{aligned}
&D \left( y-\delta \right) & \delta \leq y \\
&0 & \text{otherwise} \\
\end{aligned}
\right.
\end{aligned}
\label{policy}
\end{equation}

$y$ is a vector with the transmission capacities a TSO would like to guarantee. $D$ is a matrix which maps the lack of transmission capacity $\left( y-\delta \right) $ to suitable units. Note, that $\delta$ and $y$ are always positive. $D$ needs to be optimized together with the amount to procure. In Fig. \ref{fig:policy}, the policy for one generation unit is illustrated. If the line rating falls below $y = 1.2$ for this specific line, the output of a certain generator is set as a piece-wise linear function of the uncertainty. Also shown is the probability density function (PDF) of the line rating. Depending on the slope of the policy and the selection of $y$, the expected costs vary. An envisaged operation scheme would distribute the calculated policies to the providers of reserves. The line rating is continuously measured using sensors and distributed to the providers. In case that the line rating $y$ is not guaranteed anymore, they activate their reserves according to their policy. The operation is done in a decentralized manner. The policies could be updated based on fixed time-intervals, e.g. after market clearings, or event-based, e.g. after substantial changes to the system state have happened.

\newcommand\gauss[2]{1/(#2*sqrt(2*pi))*exp(-((x-#1)^2)/(2*#2^2))} 

\begin{figure}[!h]
\center
\begin{tikzpicture}
\begin{axis}[every axis plot post/.append style={
  mark=none,domain=0.25:1.8,samples=50,smooth}, 
axis x line=bottom, 
axis y line=left, 
enlargelimits=upper,
xlabel=DLR/NLR,
ylabel=PDF or Generator output,
ymin=-0.2,ymax=4.5,
yticklabels={,,},
legend style={at={(0.5,1.2)},
		anchor=north, legend columns=-1},
width=0.45\textwidth,
height=0.3\textwidth,
extra x ticks={1.2}, extra x tick labels={},] 
\addplot[draw=black] {\gauss{1.3}{0.1}};
\addplot[draw=red] coordinates{(0.25,2)(1.2,0)};
\addplot[draw=red] coordinates{(1.2,0)(1.8,0.01)};
\addplot[draw=black,dashed] coordinates{(1.2,0)(1.2,2)};
\legend{PDF,Affine Policy}
\end{axis}
\end{tikzpicture}
\caption{Affine policy for a generator reacting on forecast errors of a certain line rating.}
\label{fig:policy}
\end{figure}
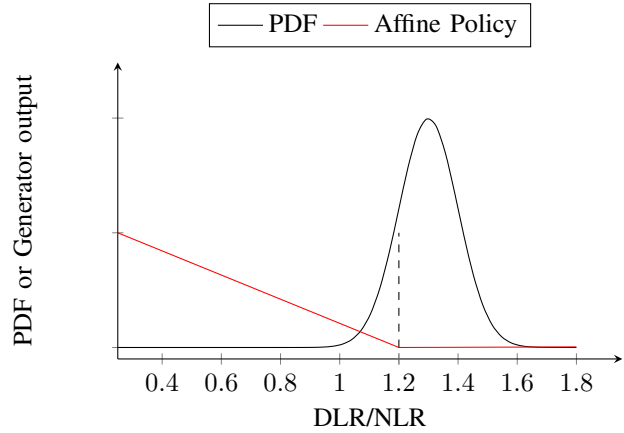

\subsection{Expected Operation}

If the current line rating $\delta$ falls below the guaranteed transmission capacity $y$ of the line, remedial actions have to be taken. In order to estimate the expected changes in the setpoints which translates in additional operational costs, the expected value of the multiplication of the policy $D$ with distribution of $\delta$ is calculated. For a normally distributed $\delta$, the expected value of this truncated distribution can directly be expressed as \cite{KotzContUniDistributions}:

\begin{equation}
\begin{aligned}
De &=D\mathbb{E}\left[ \Delta \right] = D \mathbb{E} \left[ y-\delta \left| \delta \leq y \right. \right] \Phi \left( \delta \leq y \right) \\
\mathbb{E} \left[ y - \delta \left| \delta \leq y \right. \right] &= y - \mu + \sqrt{\text{diag}\left( \Sigma \right)} \frac{\phi \left(\beta \right)}{\Phi \left( \beta \right)} \\
\beta &= \frac{y-\mu}{\sqrt{\text{diag}\left( \Sigma \right)}} \\
\end{aligned}
\label{eq:app2}
\end{equation}

$\phi\left( . \right)$ is the probability density function, and $\Phi \left( . \right)$ the cumulative distribution function of the Gaussian distribution.


\subsection{Procurement and Dispatch}

For a selected $y$ and the corresponding precalculated $e = \mathbb{E}\left[ \Delta \right]$, the problem can be formulated as:

\begin{equation}
\begin{aligned}
c_{DLR}^* =& \min_{P_{gen},\underline{\Delta},\overline{\Delta},D_{up},D_{dn}}\  c_{disp} + c_{proc} + c_{op}\\
&s.t.\ \forall \delta \in \mathcal{W}\\
\circled{A}&\left\{
\begin{aligned}
& c_{disp} \geq P_{gen}^T Q_{gen} P_{gen} + c_{gen}^T P_{gen} \\
& \underline{P}_{gen} \leq P_{gen} \leq \overline{P}_{gen} \\
& \sum_j \left( P_{gen,j} - P_{load,j}\right) = 0 \\
\end{aligned}
\right.\\
\circled{B}&\left\{
\begin{aligned}
& c_{proc} \geq c^{+} \overline{\Delta} - c^{-} \underline{\Delta} \\
& c_{op} \geq c_{op}^{+}D_{up}e + c_{op}^{-}D_{dn} e \\
& D = D_{up} - D_{dn}\\
& \sum D \delta = 0 \\
& \underline{\Delta}^{max} \leq \underline{\Delta} \leq D \left( y- \delta \right) \leq \overline{\Delta} \leq \overline{\Delta}^{max}\\
& D_{up}, D_{dn}, \overline{\Delta}, -\underline{\Delta} \geq 0 \\
\end{aligned}
\right.\\
\circled{C}&\left\{
\begin{aligned}
& P_L = H \left( P_{gen} - P_{load}\right) \\
& \left| P_{L,DLR} + H_{DLR} D \left( y- \delta \right) \right| \leq \delta\\
& \left| P_{L,NLR} + H_{NLR}D\left( y- \delta \right) \right| \leq \overline{P}_{L,NLR}\\
\end{aligned}
\right.
\end{aligned}
\label{eq:opt2}
\end{equation}


where $\mathcal{W} = \left\{z: \mu + Bz, \lVert z\rVert_2 \leq \rho, Bz \leq y - \mu \right\} $ is an uncertainty set based on a constrained ellipsoidal uncertainty set that only considers realizations where a remedial action is required.
The optimization tries again to minimize the costs for the dispatch, procurement and the operational costs of the reserves, which can be estimated due to the affine policies. As in optimization \eqref{procI}, the same constraints for the dispatch problem apply (constraints \circled{A}). The remedial actions are according to the policy \eqref{policy}, with the condition, that there is no net injection (constraints \circled{B}). The power flows on lines with dynamic line ratings have to be smaller than the actual line rating and the line limits of all remaining lines have to be satisfied after remedial actions have been taken (constraints \circled{C}). The amounts to procure are such that there is sufficient reserves available for every realization $\delta$ with the corresponding $D$. The robust optimization problem Eq. \eqref{eq:opt2} can be solved efficiently as shown in \cite{BenTal99}.



\clearpage{}
\section{Case Study}
\begin{figure}[!h]
\centering
\includegraphics[width=.48\textwidth]{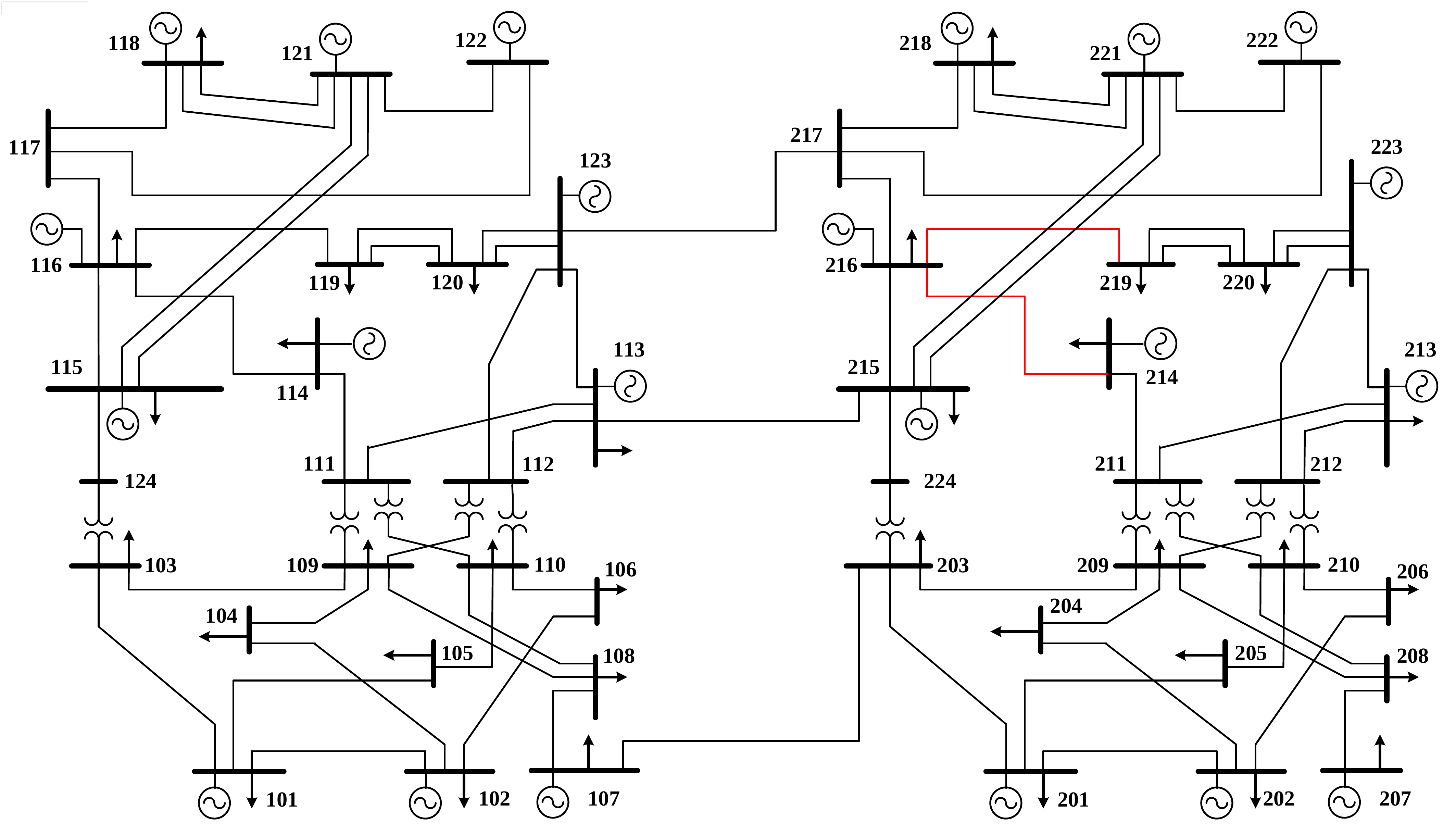}
\caption{RTS96 grid with 2 areas \cite{RTS96}.}
\label{RTSFIG}
\end{figure}

We apply the methods presented above in a case study and compare the results and the influence of certain parameters.
As test grid we use the IEEE RTS-96 2 zones test system shown in Fig. \ref{RTSFIG}. The lines marked red, 214-216 and 216-219, are selected to have dynamic line ratings and their nominal line ratings are set to $1 p.u. = 250\textrm{MW}$.

Exemplarily, we pick a scenario, where the point forecast of the line ratings are  $\mu = 1.5$p.u. The forecast error for the day-ahead point forecast is $\sigma_{24} = 0.2$p.u. and for the 3h-ahead forecast $\sigma_{3} = 0.1$p.u. Without loss of generality, the off-diagonal elements of the covariance can be set to zero, i.e. the line ratings are not correlated. It should be noted that in general forecast errors are a function of the point forecast itself, but this is neglected for reasons of simplicity \cite{MB2}.

\subsection{Increasing social welfare}
\begin{figure}
\centering
\begin{minipage}{.4\textwidth}
  \centering
  \includegraphics[width=.8\linewidth]{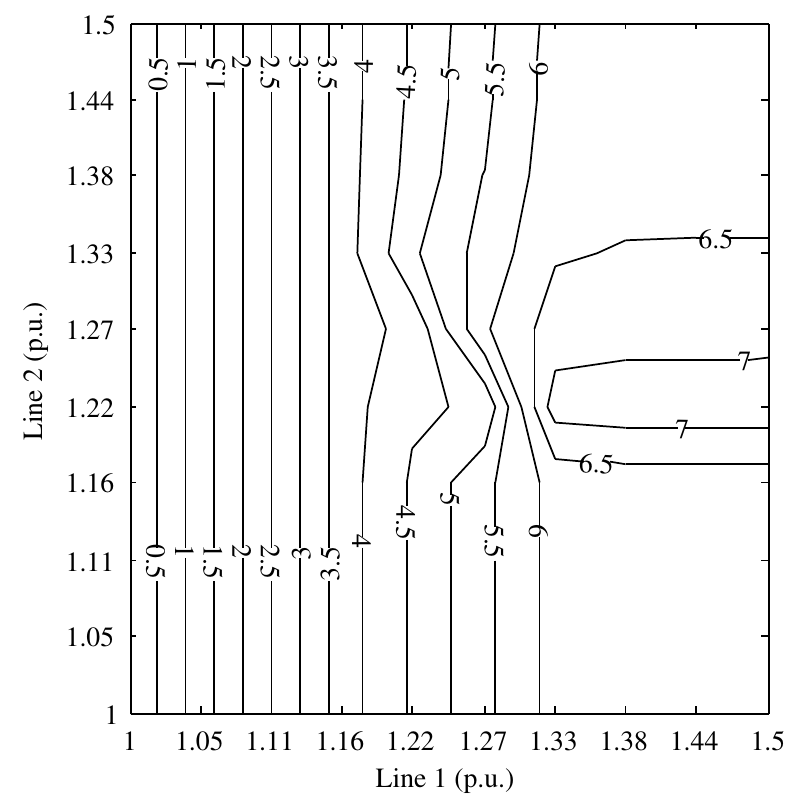}
  \captionof{figure}{Cost savings in percent of total operational costs for 3h-ahead forecast.}
  \label{fig:c3a}
\end{minipage}%
\hspace{1cm}
\begin{minipage}{.4\textwidth}
  \centering
  \includegraphics[width=.8\linewidth]{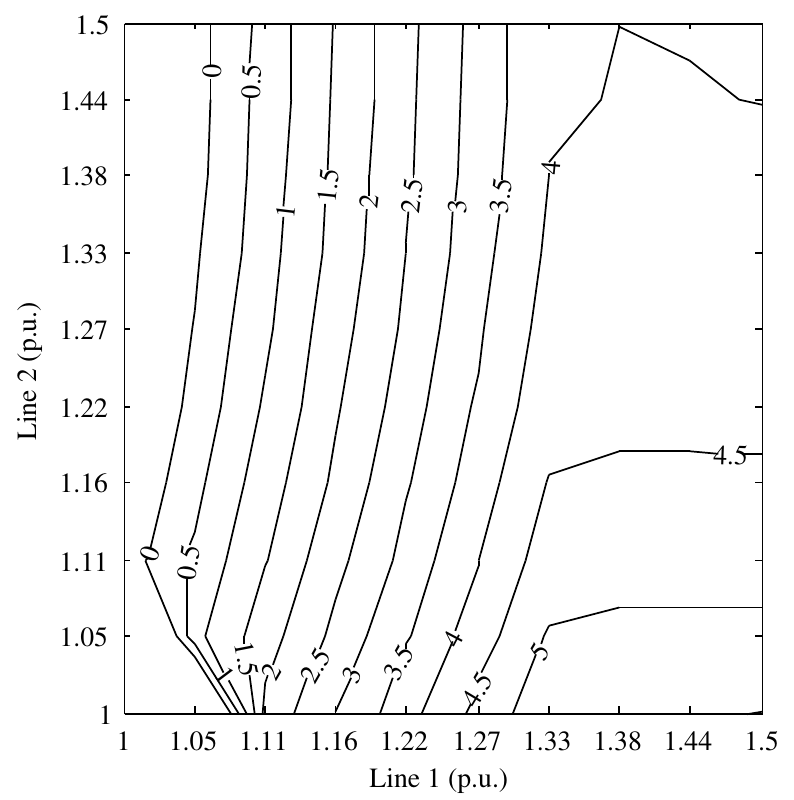}
  \captionof{figure}{Cost savings in percent of total operational costs for 24h-ahead forecast.}
  \label{fig:c24a}
\end{minipage}
\end{figure}

Additional transfer capacity enables a more flexible operation of the transmission system and frees capacity, e.g. for cross-border balancing and trading.
Therefore, reduced overall costs are expected as less expensive generation units can be dispatched which were limited before due to congestions.
In this case study, the goal is to demonstrate the reduction of total costs for the operation of the power system. We consider a 3h-ahead and a 24h-ahead forecast with the same expected value but with different uncertainties as stated above. We provide the results only for approach II, as the expected costs for operating the DLR system are directly incorporated in the optimization. Using a Monte Carlo based simulation, similar results can be produced for approach I. Figs. \ref{fig:c3a} and \ref{fig:c24a} show contour plots with the percentual reduction in costs for different maximum allowed flows on the two lines with DLR. The following is observed: first, when the forecast uncertainty is reduced, i.e. moving from 24h-ahead to 3h- ahead forecast, the cost reductions are higher. Further, the plots show that the influence of Line 2 on the operational costs are smaller compared to Line 1. Therefore should focus only on Line 1 for the given conditions. For the given costs for procurement and operation as well as the forecast, the costs are reduced in all cases. However, for different cost structures, or if different lines are equipped with DLR and different forecast quality, the results may change. Future work should therefore investigate limitations on economical operation of DLR.

%
%

\begin{table*}[!hbp]
\caption{Comparison of approach I and II for given forecast and different parameters.}
\begin{tabular}[C]{|l|l|l|l|l|l|l|l|}
  \hline
   & status quo& II-3h & II-24h & I-3h ($\alpha = 0)$ & I-3h ($\alpha = 0.1)$ & I-24h ($\alpha = 0)$ & I-24h ($\alpha = 0.1)$ \\\hline
  Dispatch Cost & 108760 &100770&100770&100720&100770&100720&100770 \\\hline
  Procured Amount & 0 & 110.35&330.89&10.45&12.05&68.75&134.28\\\hline
  Procurement Costs & 0 & 1103.5&3308.9&104.5&120.5&687.5&1342.8\\\hline
  Mean operational costs& 0 & 33.26&92.47&22.17&0.5&673.92&44.11\\\hline
  Total Cost & 108760 & 101906&104171&100846&100891&102081&102156\\\hline
 Cost Saving (in \%) & 0 & 6.3 & 4.2 & 7.3 & 7.2 & 6.1 & 6\\
  \hline
\end{tabular}
\label{tab}
\end{table*}

Tab. \ref{tab} compares approach I and II for the given forecast and $y$ is selected such that the cost reductions are the highest. In this snapshot, we compare both methods with forecast lead times. Additionally, method I is tested for two values of the weighting factor $\alpha$ (see Eq. \eqref{procI}).
It can be seen that for the given snapshot the cost savings are all positive, i.e. the total costs for operating the power system are reduced. The procured amounts of reserves are varying a lot, but they mainly depend on the selection of the risk level $\gamma$ and the selected line. The results in approach I further depend on $\alpha$ and in approach II on $y$. We observe two tendencies: first, increasing the $\alpha$ in approach I increases the procured amounts, but as the worst-case operational costs are considered, the average operational costs determined by a Monte Carlo simulation are lower. Second, the forecast period plays a crucial role, i.e. the proposed approaches performs better for short-term weather forecasts and thus in an intra-day operation scheme. It could be considered as if the TSO would release additional transmission capacity to the intra-day market if it is available. Day-ahead markets would be cleared respecting the nominal line ratings. Lastly, the costs are higher for procurement and operation of the DLR for approach II. As approach II works in a distributed manner, it might be used for small deviations and only with few lines equipped with DLR.

\subsection{Procured Amounts}
\begin{figure}[!h]
\centering
\includegraphics[width=.4\textwidth]{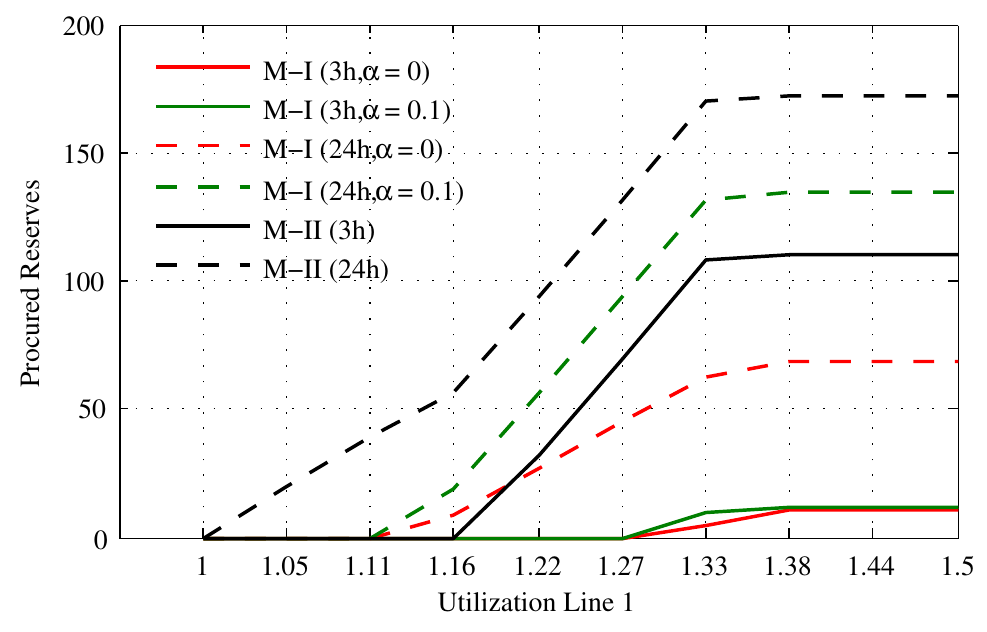}
\caption{Procured amounts for different upper flow limits on Line 1.}
\label{fig:proc}
\end{figure}

Fig. \ref{fig:proc} shows the sum of procured amounts for different upper flow limits on line 1 for the given weather forecasts. We apply for both approaches the 3h forecast and 24h forecast case. In the case of approach I, the $\alpha$ is selected to be 0 and 0.1, i.e. the worst-case operational costs are considered. One can observe that the required reserves increase with increasing line utilization and larger forecast lead-times. In almost all cases, the line could be utilized more than the nominal line rating without the need to procure reserves. This however depends strongly on the given forecast. Further, it can be seen, that the amounts reduce if the forecast lead-time and thus the error is reduced. Comparing the approaches I and II, one can see that approach I requires less reserves as it is not subject to a specific policy. If $\alpha$ is selected as $0.1$, the amounts are increasing slightly. However, one should consider, that in such a case the operational costs are likely to be smaller. It should be noted, that aside from the procured amounts also the location where the reserves are located are crucial, which is considered via the grid constraints in the optimizations.

\subsection{Influence of Forecast Uncertainty}
As the distribution of the forecast error is in general a function of the point forecast itself, we investigate the costs for the operation of the reserves as well as the reduction in the overall costs for different forecasts and different forecast errors. We select three forecasts, $\mu = \left\{1.25, 1.5, 1.75 \right\}p.u.$, and vary the diagonal elements of the covariance $\Sigma$ between $0.01^2 p.u.$ and $0.3^2 p.u.$. The off-diagonal elements are zero.

\begin{figure}
\centering
\begin{minipage}{.4\textwidth}
  \centering
  \includegraphics[width=.8\linewidth]{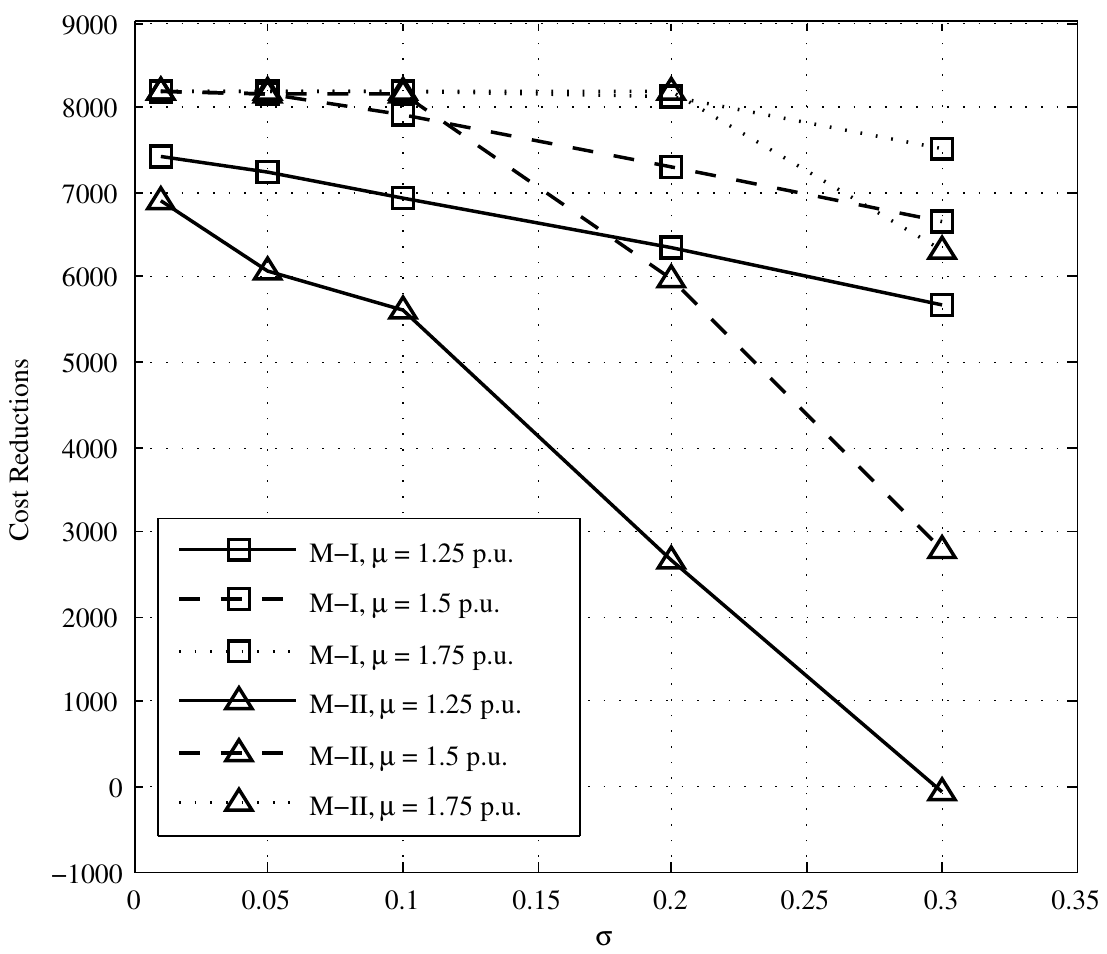}
  \captionof{figure}{Cost Reductions for different $\mu$ and $\sigma$ per timestep.}
  \label{fig:musigma1}
\end{minipage}%
\hspace{1cm}
\begin{minipage}{.4\textwidth}
  \centering
  \includegraphics[width=.8\linewidth]{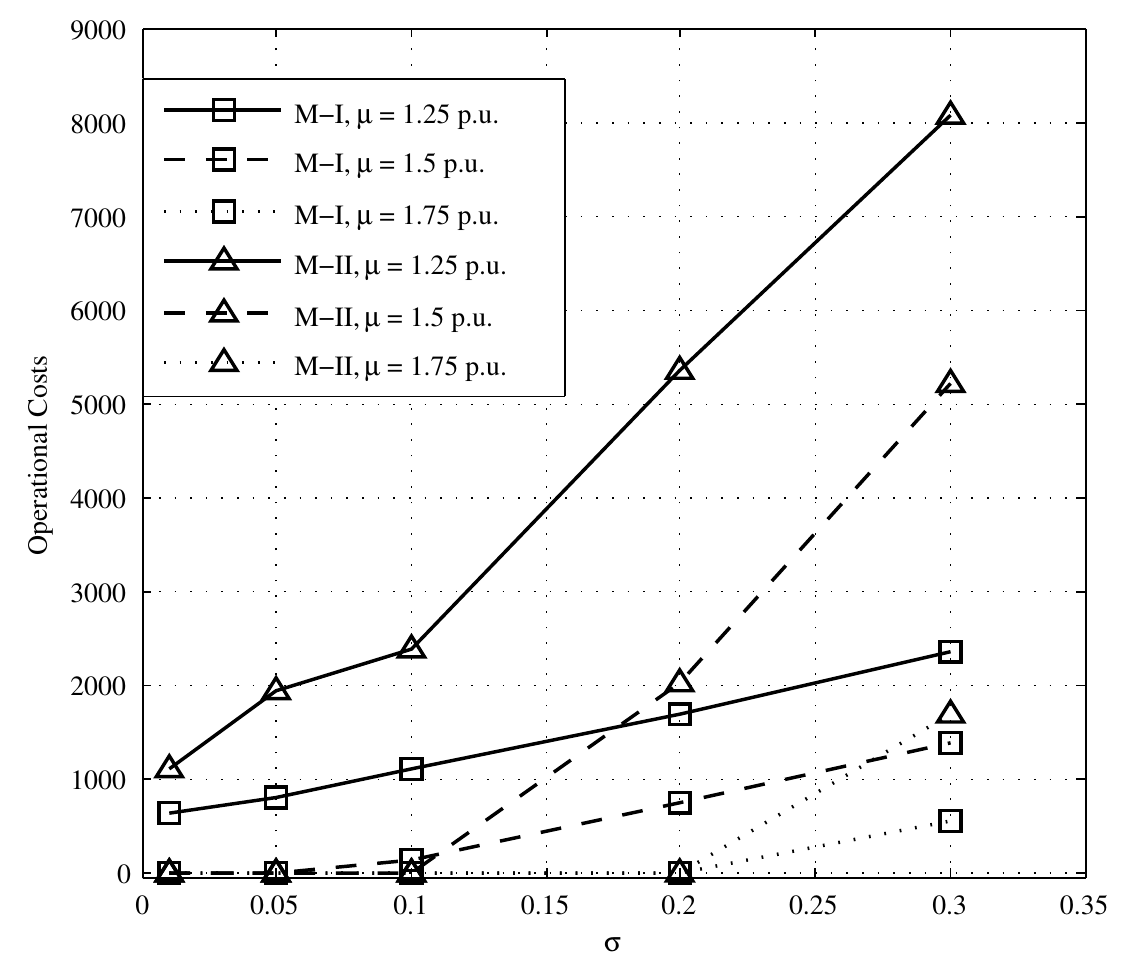}
  \captionof{figure}{Operational costs for different $\mu$ and $\sigma$ per timestep.}
  \label{fig:musigma2}
\end{minipage}
\end{figure}

%

Fig. \ref{fig:musigma1} shows the average cost reductions in total operation costs including costs for DLR operation for both approaches and the three different forecasts. Similarly, Fig. \ref{fig:musigma2} shows the average costs arising from remedial actions due to forecast errors. The average costs are determined by applying the methods to 1000 scenarios. For comparison, the total dispatch costs without DLR for the considered time period is 108'760.\\

The results are as expected: the operational costs increase with increasing uncertainty. Accordingly, the cost reductions are lower. It is therefore desireable to operate on a short-term basis in order to keep the forecast errors small. A higher forecast leads to a general higher utilization and therefore increasing cost reductions and lower operational costs. Approach I has lower operational costs and higher cost reductions and is thus the better option from an economic point of view.

\subsection{Discussion}
Approach I has the advantage that the procurement is easy and computationally fast. The operation, for a given procured amount, is cost-efficient, as an online optimization is used. The drawbacks are that the operation needs centralized operation and communication infrastructure as the re-dispatch setpoints have to be communicated to the units. Further, the online optimization might complicate the time-sensitive operation of the units. Finally, the procurement costs has only limited information on the expected operational cost and thus might lead to higher operational costs. Approach II enables a distributed operation as the changes of the units' setpoints are directly related to the thermal line rating measurements. The policies can be precalculated and the expected operational costs are integrated in the optimization but the computation is harder than for approach I and might lead to an inefficient operation when the forecast errors are substantially large.
Both approaches can enable cost savings, but the centralized approach I gives higher savings in general. Therefore, the topic of future research should address a combination of both methods, e.g. considering approach II as fallback option when communication to the coordinating entity is disrupted.\\

Aside from the thermal limits of the lines, also other factors neglected in this paper, should be considered in a implementation, such as voltage limitations, the N-1 security criterion or operational limits of other equipments, e.g. transformers. For a successful implementation of DLR, also a real-time monitoring and reliable models of the thermal state of transmission lines are crucial \cite{Sattinger}.

\section{Conclusion}
In this paper, two approaches that enable a robust operation of power systems with transmission lines using dynamic line ratings have been presented. Both approaches calculate the amount and location of power that can be re-dispatched. The calculations are based on a forecast of the meteorological conditions and possible forecast errors given by an uncertainty set. The first approach relies on centralized optimization whereas the second approach works in a decentralized manner using affine policies. Case studies compared the two approaches and demonstrated the reductions on total operational costs as well as the necessity to operate such a system on a short-term basis. Further, variation of statistical parameters have shown, that precise forecasts of meteorological conditions are crucial.\\
Future research should investigate DLR in combination with controllable devices such as HVDC or FACTS as well as topology changes, e.g. line switching. The required number of measurement devices and their placement in the grid is another open question. Further, the incorporation of non-normally distributed forecast errors in given methods as well as their application to realistic grid setups using historic meteorological data should be simulated in order to quantify the cost reductions.

\footnotesize
\section*{Acknowledgement}
This research was carried out within the project Balancing Power in the European System (BPES). Financial support by the Swiss Federal Office of Energy (SFOE) and Swissgrid is gratefully acknowledged.

\footnotesize
\bibliographystyle{IEEEtran}
\bibliography{Bibliography}
\end{document}